\documentclass[11pt]{article}

\usepackage{epsf}
\usepackage{psfrag}
\usepackage{amsmath}
\usepackage{amssymb}
\usepackage{lscape}
\usepackage{latexsym}
\usepackage[usenames]{color}
\usepackage[mathcal]{eucal}
\usepackage{stmaryrd}
\usepackage{textcomp}
\usepackage{setspace} 
\usepackage{graphicx} 
\usepackage{pgfplots} 
\usepackage{epsfig}   
\usepackage{epstopdf} 
\usepackage{appendix}
\usepackage{afterpage}
\usepackage{float} 
\usepackage[misc,geometry]{ifsym} 

\newtheorem{theorem}{Theorem}[section]

\newtheorem{example}{Example}[section]

\newcommand{\ds}{\displaystyle}
\newcommand{\dZ}{{\cal Z \kern -0.7em Z}}
\newcommand{\dC}{{\rm\hbox{C \kern-0.8em\raise0.2ex\hbox{\vrule height5.4pt width0.7pt}}}}
\newcommand{\dQ}{{\rm\hbox{Q \kern-0.85em\raise0.25ex\hbox{\vrule height5.4pt width0.7pt}}}}

\newcommand{\proofbox}{\hspace{\fill}{$\Box$}}

\newenvironment{proof}{Proof.}{\proofbox}

\newcommand\old[1]{}
\newcommand{\beqa}{\begin{eqnarray*}}
\newcommand{\eeqa}{\end{eqnarray*}}

\title{\bf \Large {A Multiobjective Mathematical Model of Reverse Logistics  for Inventory Management with Environmental Impacts: An  Application in Industry}}

\author{M. Forkan\thanks{Department of Mathematics,
    University of Chittagong, Chittagong-4331, Bangladesh,
   {\tt email: forkan.math@cu.ac.bd},}
\and
M. M. Rizvi \Letter \thanks{Department of Mathematics,
	University of Chittagong, Chittagong-4331, Bangladesh, and University of South Australia-STEM, Adelaide {\tt email: rizmm001@mymail.unisa.edu.au} (corresponding author),}
\and
M. A. M. Chowdhury\thanks{JNIRCMPS,
	University of Chittagong, Chittagong-4331, Bangladesh, {\tt email:
      mamansur@cu.ac.bd}}}
\begin{document}
%
\maketitle
\textbf{Abstract} We propose new mathematical models of inventory management in a reverse logistics system. The proposed models extend the model introduced by Nahmias and Rivera with the assumption that the demand for newly produced and repaired (remanufacturing) items are not the same. We derive two mathematical models and formulate unconstrained and constrained optimization problems to optimize the holding cost. We also introduce the solution procedures of the proposed problems. The exactness of the proposed solutions has been tested by numerical experiments. Nowadays, it is an essential commitment for industries to reduce greenhouse gas (GHG) emissions as well as energy consumption during the production and remanufacturing processes. This paper also extends along this line of research, and therewith develops a three-objective mathematical model and provides an algorithm to obtain the Pareto solution.


\textbf{AMS} subject classifications. 90B05, 90C25, 90C29, 90C30.

\section{Introduction}\label{sec1}
Reverse logistics (RL) has been defined as a term that refer to the role of logistics in product returns, a recycling, the materials substitution, a reuse of materials, a waste disposal, and the refurbishing \cite{Govindan 2017}. Moreover, inventory management in reverse logistics, which incorporates joint manufacturing and remanufacturing options, has received increasing attention in recent years. However, fast developments in technology and mass appearance of new industrial products, which are coming to the market, have resulted in an increasing number of idle products and caused growing environmental problems worldwide. Therefore, increasing ecological concerns, end user awareness, economic considerations, and legislation, related to waste disposal, encourage manufacturers to take back products after customer have used them. Recently, growing interest and realizations in the reverse logistics processes, such as the recovery of the returned products, have become one of the ways, in which businesses endeavor to retain and increase competitiveness in the global market.

 \par 
A number of publications with various models have appeared in the literature aimed to optimize holding cost in the process of reverse logistics system. McNall \cite{McNall1966} and Schrady \cite{Schrady1967} were the first to address the inventory problem for repairable (recoverable) items.  Schrady \cite{Schrady1967} first explored a deterministic reverse logistic Economic Order Quantity (EOQ) model for repairable items with multiple repair cycles and one production cycle. The model of Schrady \cite{Schrady1967} was extended by Nahmias and Rivera    \cite{Nahmias79} with inclusion of the case of finite repair rate. Richter (\cite{Richter94}-\cite{Richter96b}) proposed an EOQ model with waste disposal and looked over the optimal figure of production and remanufacturing batches, depending on the rate of return. Richter \cite{Richter97}, Richter and Dobos \cite{Richter99, Dobos2000} investigated whether a policy of either total waste disposal or no waste disposal is optimal. Teunter \cite{Teunter2001} considered multiple productions and remanufacturing cycles and generalized the results from Schrady \cite{Schrady1967}. Dobos and Richter \cite{Dobos2003} developed a production/recycling setup with constant demand that is satisfied by non-instantaneous production and recycling with a single repair and a single production batch in an interval of time. Later on, Dobos and Richter \cite{Dobos2004} generalized their earlier model \cite{Dobos2003} by considering multiple refurbish/repair and production batches in a time interval. Along the same line of study, Dobos and Richter \cite{Dobos2006} further extended the model and assumed that the quality of collected used/returned items is not always suitable for further recycling. Later on, Jaber and El Saadany \cite{Saadney09, Saadney10} extended the work of Richter \cite{Richter96a,Richter96b} by assuming that the remanufactured items are considered by the customers to be of lower quality than the new ones. Alamri \cite{Alamri2011} put forward a general reverse logistics inventory model for deteriorating items by considering the acceptable returned quantity as a decision variable.  Singh and Saxena \cite{Singh2012} proposed a reverse logistics inventory model allowing for back-orders. Hasanov et al. \cite{Hasanov2012} extended the work of Jaber and El Saadany \cite{Saadney09} by assuming that unfulfilled demand of remanufactured and produced items is either fully or partially backordered. El Saadany et al. \cite{Saadney13} discussed an inventory model with the question as to how many times a product can be remanufactured. Singh and Sharma \cite{Singh2013b} explored an integrated model with variable production and demand rates under inflation. Later, Singh and Sharma \cite{Singh2016} established a production reliability model for deteriorating products with random demand and inflation. Recently, Bazan et al.  \cite{Bazan2015, Bazan16} presented a  mathematical inventory models for reverse logistics with environmental perspectives.

\par 
In the existing literature, most of the research articles are developed with the assumption that the produced and recovered items are not of different quality. In many practical business situations this hypothesis is not adequate, as the repaired (remanufactured) items are considered of secondary quality by the customers, for instance, wheel tyre and computer, etc. Therefore, in this study, it is assumed that newly produced and remanufactured items are different in quality. This paper is an extension of the work of Nahmias and Rivera \cite{Nahmias79} for the case of finite repair rate, production and remanufacturing. Moreover, the solution approach of our proposed model enriched and propagated to the case where constraints of confined storage space in the repair and supply depots are imposed. Note that the portions of products to be procured, repaired, and disposed of in per time unit are fixed, so that, the unit cost of these does not have any impact on optimizing the reverse logistics model. Therefore, our proposed model is associated only with setup and holding costs of items of supply and repair depots. Numerical experiments are provided to illustrate the proposed models. The behaviors of the total average holding cost functions for different stock-out cases are presented with respective tables and graphs.

 \par 
In addition, we developed a multiobjective mathematical model of the reverse logistics
system that satisfies not only holding cost objective but also environmental requirements.
The model considers two environment perspectives; one is to minimize the greenhouse gas
emission during the production process, and the other is to optimize energy used in the production and remanufacturing processes. Since the problem has complex nature; the effective methodology and efficient
algorithms are needed to approximate the solutions on the front. So far we know that very few researchers have been devoted to solving multiobjective reverse logistics model \cite{Bazan2015, Bazan16}. One of the reasons for not having enough literature about the solution approach is that the feasible set is not convex and might even be disconnected due to the conflicting nature of the multiple objectives. This poses difficulties for any technique to approximate the Pareto solutions. We initiate a well-known scalarization approach and algorithms \cite{BurKayRiz2017} to solve the proposed three-objective reverse logistics model. Extensive numerical experiments are conducted to approximate the non-dominated solution of the multiobjective model and the Pareto solutions have been demonstrated. \par

 In this paper, we have also integrated the applications of our models  in the tyre industry. In the computational experiments, we consider three numerical examples that are presented in such a way so that they can align with the industries where the demands of new and repaired items are different (for example, tyre and IT industries, etc.). At the same time, we would like to test our proposed models' capabilities and analyze the obtained solutions under different parameter settings.  We are aware that in the tyre industry the requirements of new and repaired tyres are always different. Customers have different choices to use new and repaired tyres. Some customers prefer the tyre's quality, and they do not compromise with the safest driving conditions, so they can make a choice for a new tyre, and others can choose the repaired tyre as the repaired tyre is cheaper and fulfill buyers' requirements.

The rest of this paper comprises the following sections. In Section \ref{sec2}, the new mathematical reverse logistics models for unconstrained and constrained optimization problems and their solutions are illustrated. The extension of these models in the context of multiobjective case is
described in Section~\ref{sec3} and in Section~\ref{sec3a}, a solution approach for three objectives problem is proposed. In Section~\ref{sec4}, numerical experiments are conducted and we provide the results and discussions based on numerical experiments. The last section presents the conclusion of the paper.
\section{Formulations and Solutions of the Proposed Models }\label{sec2}
\begin{figure}[h]
	\centering 
	\includegraphics[width=1.10\textwidth]{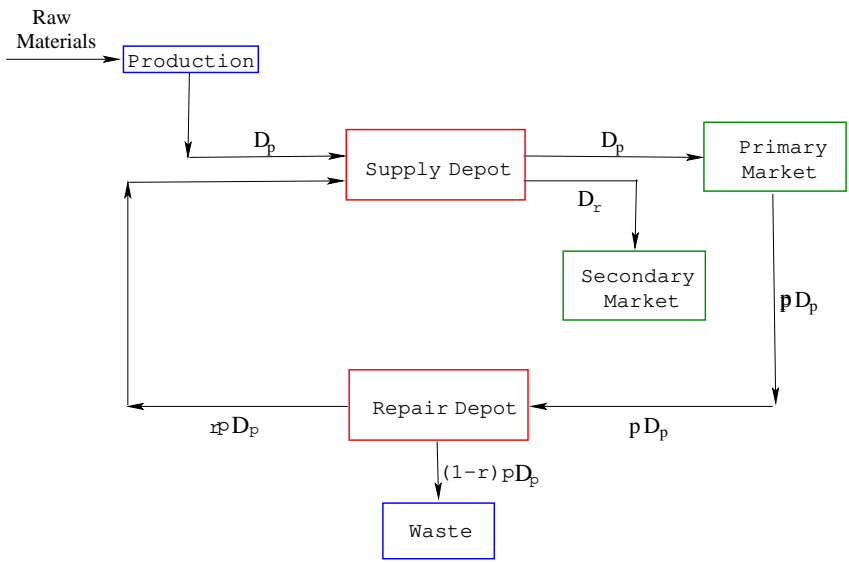}
	\caption{Material flow for a new procurement and repair system.}
	\label{FlowModel}
\end{figure}

We consider two type of depots in the proposed model. The first depot is the supply depot
which stores the new procurement and repaired items. The user's demand of primary and secondary markets can be satisfied from this supply depot. When the items are repaired,
first they are stored at the repair depot (second depot) and subsequently shipped in batches
to the supply depot. We assume that shortage is not allowed in this stock point and the lead time is zero.The demand of the new and repaired items are fixed in time but may not be equal, as the repaired items sell at lower price in the secondary market. Let $D_p$ and $D_r$ are the demands for new and repaired items respectively in per unit time. Note that the used items are sent back from user to the overhaul and then repair depot with a constant rate. The material flow of the model is
depicted in Figure~\ref{FlowModel}. \par  

It is assumed that the procurement and repair batch sizes are $Q_p$ and  $Q_r$, respectively over the time cycle $T$.  The following parameters have been used in the model.  The collection percentage of available returns of used items is $p$ $(0<p<1)$, the recovery rate $r$, waste disposal rate $1-r$, $\lambda$ be the repair rate per unit time with $\lambda > D_r$, fixed procurement cost $A_p$, fixed repair batch induction cost $A_r$, holding costs at supply and repair depots are $h_1$ and $h_2$ per unit per unit of time, respectively. \par 

\begin{figure}[h]
	\centering 
	\includegraphics[width=.70\textwidth]{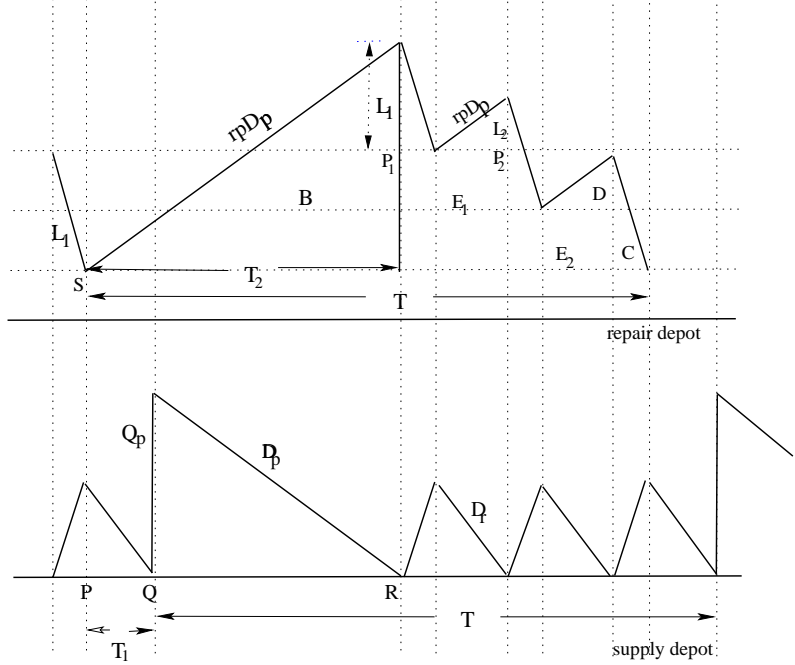}
	\caption{The behaviour of inventory for  produced, collected used and repaired items over interval $T$.}
	\label{inventoryModel}
\end{figure}

In the proposed models, we consider $n$ cycles repair and single-cycle production over the time $T$.
Therefore, the relations between in and outflows of the stocking points in a procurement and $n$
repair cycles are \[ Q_p+n.Q_r=(D_p+D_r).T\;\;\text{and}\] \[n.Q_r=r.p.D_p.T\]
The behavior of inventory for produced, collected used and repaired items over the time
interval $T$ is illustrated in Figure~\ref{inventoryModel}. According to the Figure~\ref{inventoryModel}, the inventory of the supply depot always decrease at rates $D_p$ and $D_r$ for new and repaired items, respectively. The inventory of the repair depot increase at rate rpDp, and pDp is the portion of the demand which is returned to the system for either repair or dumping. The portion of demand $(1-p)D_p$ never return from primary market to the system, so we assume that it does not have any influence in the system. The returned items are usually of varying quality. It is considered that a returned item with a quality less than the acceptance quality is rejected. In the process, the amount of the failed/waste item is $(1-r)pD_p$. The decrease of 
repaired items during repair is $\ds L_1=\left(\frac{Q_r}{\lambda}\right)(\lambda-rpD_p)$. It is mentioned here that, we offen introduce  various notations in the text to derive the simplest expression for holding cost function. Let,  $\ds C_1=\left(1-\frac{rpD_p}{\lambda}\right)$, where $\lambda >rpD_p$, this follows that $L_1=C_1Q_r$. The repaired item falls below the point $S$, the repair is suspended, and the inventory in the repair depot grows at a rate of $rpD_p$ for a time $T_2$ (see, Figure~\ref{inventoryModel}). We obtain $L_1-T_1D_r=0$, thus,
\begin{equation*}
T_1=\frac{C_1Q_r}{ D_r}.
\end{equation*}
At the time $T_1$ a new procurement $Q_p$ is received in the supply depot and is used to meet
demand whereas return items accumulated at the repair depot. At the time $T_2$ when
the new procurement is all distributed, repaired items are initiated at the supply depot.
Therefore, we have
\begin{equation*}
T_2=PQ+QR=T_1+\frac{Q_p}{D_p},
\end{equation*}
it follows that,
\begin{equation*}
T_2=\frac{C_1Q_r}{ D_r}+\frac{Q_p}{D_p}.
\end{equation*} 
The cycle length $T$ be the time between two new procurements which is also the same as the time between two successive suspensions of repair.  The items fail to repair for each successive induction at a constant rate $(1-r)pD_p$ for a period of time $P_1P_2$, that is, $\ds T_1+\frac{Q_r}{\lambda}$. In this period the net loss is $L_1-L_2$ (see, Figure-2) which is $C_1Q_r-T_1rpD_p$ . The total loss of returned item from the system due to item failed to qualify the quality test is $(n-1)(C_1Q_r-T_1rpD_p)$. It can be seen from Figure~\ref{inventoryModel} that the following relation holds as
\begin{equation*}
(n-1)(C_1Q_r-T_1rpD_p)=rpD_pT_2-Q_rC_1,
\end{equation*}
this gives 
\begin{equation}\label{n_value}
n=\frac{C_2Q_p}{Q_r},
\end{equation}
where $\ds C_2=\frac{rp}{C_1\left(1-\frac{rpD_p}{D_r}\right)}$. The total number of units giving up the supply depot in a cycle $T$ is exactly $T(D_p+D_r)$ which must be equal to the total number of units ingoing into the supply depot in the same time, which is $Q_p+nQ_r$, therefore, 
\begin{equation}\label{T11_value}
T=\frac{Q_p+nQ_r}{D_p+D_r}.
\end{equation}
From (\ref{n_value}) and (\ref{T11_value}) we have, $T=C_3Q_p$, where $\ds C_3=\frac{1+C_2}{D_p+D_r}$. \\ [1.5ex]
Now we evaluate the area of the inventory curve of the supply depot during a cycle $T$.
Let $A_1$ be the area bounded by the inventory curve of the supply depot over the time T. Thus, from Figure~\ref{inventoryModel} we get  
\begin{equation}\label{supplyArea}
A_1=\frac{Q_p^2}{2D_p}+\frac{n}{2}L_1\left[T_1+\frac{Q_r}{\lambda}\right],
\end{equation}
after substituting $L_1$, $T_1$ and $n$ in (\ref{supplyArea}), gives the area of the supply depot as
\begin{equation*}\label{A1}
A_1=\frac{Q_p^2}{2D_p}+\frac{C_1C_2Q_pQ_r}{2}\left[\frac{C_1}{D_r}+\frac{1}{\lambda}\right].
\end{equation*}
We assume the inventory level in the repair depot is zero when the repair process is suspended. Now we compute the area formed by the inventory curve in the repair
depot to find the associate holding costs. Since we are adopting that starting and ending
values on this curve are zero, it ensues that the bounded area in the repair depot can be
parted into triangles and rectangles as displayed in Figure~\ref{inventoryModel}.
 \par 
\noindent The area of triangle $B$ is $\ds \frac{1}{2}rpD_pT_2^2$, by substituting $T_2$, we have
\begin{equation*}\label{B}
B=\frac{1}{2}rpD_p\left(\frac{C_1Q_r}{ D_r}+\frac{Q_p}{D_p}\right)^2.
\end{equation*} 
Since $n$ cycles of repair inducted in a interval of time $T$, therefore, there are $n$ triangles $C$ with the total area $\ds \frac{n}{2}\frac{Q_r}{\lambda}L_1$,
and denoted by,
\begin{equation*}\label{C}
C^{\prime}=\frac{1}{2\lambda}C_1C_2Q_pQ_r.
\end{equation*}
We also have $n-1$ triangles $D$ with the total area $\ds \frac{n-1}{2}rpD_pT_1^2$,
and denoted by,
\begin{equation*}\label{D}
D^{\prime}=\frac{1}{2}rpD_p\left(\frac{C_2Q_p}{Q_r}-1\right)\left(\frac{C_1Q_r}{ D_r}\right)^2.
\end{equation*}
Rectangle $E_1$ has an area $\ds \left(rpD_pT_2-L_1\right)\left(\frac{Q_r}{\lambda}+T_1\right)$,
this follows,
\begin{equation*}
E_1=Q_r\left(\frac{1}{\lambda}+\frac{C_1}{D_r}\right)\left(\frac{rpD_pC_1Q_r}{D_r}+rpQ_p-C_1Q_r\right).
\end{equation*}
Rectangle $E_2$ has an area $ \ds \left(L_1-rpD_pT_1\right)\left(\frac{Q_r}{\lambda}+T_1\right)$,
thus,
\begin{equation*}
E_2=Q_r^2\left(\frac{1}{\lambda}+\frac{C_1}{D_r}\right)\left(C_1-\frac{rpD_pC_1}{D_r}\right).
\end{equation*}
Total area of the repair depot is
\begin{equation*}
A_2=B+C^{\prime}+D^{\prime}+E_1+E_2.
\end{equation*}
It is now pursues that the total average cost for set-up and holding incurred in one cycle, say $\ds f(Q_p,Q_r)$, is given by
\begin{equation*}
f(Q_p,Q_r)=A_p+nA_r+h_1A_1+h_2A_2.
\end{equation*}
Now we derive an expression of the average holding cost per unit per unit time is determined
by first finding the total cost incurred in a cycle and then dividing by the cycle
length. Hence, the average holding cost per unit per unit time, say $f_1(Q_p,Q_r)$, is created
by taking $f(Q_p,Q_r)/T$. Therefore,
\begin{equation}\label{costfunction}
f_1(Q_p,Q_r)=\frac{1}{C_3Q_p}(A_p+nA_r+h_1A_1+h_2A_2).
\end{equation}
\begin{theorem}
	$(Q_p^*,Q_r^*)$ is the global minimum of (\ref{costfunction}) occurs at
	\begin{equation}\label{partialQp}
	Q_p^*=\sqrt{\frac{2A_pD_p}{h_1+h_2pr}},
	\end{equation}
	and 
	\begin{equation}\label{partialQr}
	Q_r^*=\sqrt{\frac{2\lambda C_2A_rD_r}{C_1C_2D_r(h_1+h_2)+2D_rh_2pr+\lambda C_1\left(C_1C_2h_1+4h_2pr+\frac{C_1C_2D_ph_2pr}{D_r}\right)}}.
	\end{equation}	
\end{theorem}
\begin{proof}
Let us now take the partial derivatives of (\ref{costfunction}) with respect to $Q_p$ and $Q_r$, and then one can find the optimal procurement and repair batches by making $\ds \frac{\partial f_1}{\partial Q_p}=0$ and 
$\ds \frac{\partial f_1}{\partial Q_r}=~0$, these give \eqref{partialQp} and \eqref{partialQr}.
\end{proof}

Let us now formulate a mathematical model that minimize the average cost per unit per unit time subject to the constraints. These constraints include available floor space for the repair and supply depots. The problem needs to fulfil these requirements. Suppose $p_1$ and $p_2$ be the amount of square feet required to each of the item in supply and repair depots, respectively. We also assume that the availability of maximum number of floor space for supply and repair depots are $k_1$ and $k_2$, respectively. The highest level of supply inventory is $Q_p$, thus, the useable floor space is $p_1Q_p$ which is the less than equal of $k_1$. Similarly, the maximum level of inventory in the repair depot is the $\ds \left(\frac{C_1Q_r}{D_r}+\frac{Q_p}{D_p}\right)rpD_p$ is less than equal $k_2$. Thus, the proposed mathematical model of the problem with constraints is as follows.
\begin{equation}\label{mathemodel}
\begin{array}{rl} \min & \ f_1(Q_p,Q_r)\\
\mbox{subject to the constraints} \\
& \ds p_1Q_p \leq k_1, \\ 
& \ds p_2 \left(\frac{C_1Q_r}{D_r}+\frac{Q_p}{D_p}\right)rpD_p\leq k_2. 
\end{array}
\end{equation}

The problem is convex, therefore, sufficiency of Karush Kuhn-Tucker (KKT) conditions holds and we do not need any further conditions \cite{Miettinen1999}. Note that the objective function $f_1(Q_p,Q_r)$ and constraints are continuously differentiable.  Therefore, the KKT conditions of the Problem (\ref{mathemodel}) is presented as below.
Suppose that there exist multipliers $\lambda_1 \geq 0$ and  $\lambda_2 \geq 0$ such that 
\begin{equation}\label{KKT1}
\nabla f_1(Q_p,Q_r) + \lambda_1 \nabla (p_1Q_p -k_1) +\lambda_2 \nabla \left(p_2 \left(\frac{C_1Q_r}{D_r}+\frac{Q_p}{D_p}\right)rpD_p - k_2\right)=0,
\end{equation}
\begin{equation}\label{KKT2}
\lambda_1(p_1Q_p -k_1)=0,
\end{equation}
\begin{equation}\label{KKT3}
\lambda_2\left(p_2 \left(\frac{C_1Q_r}{D_r}+\frac{Q_p}{D_p}\right)rpD_p - k_2\right)=0.
\end{equation}
Taking partial derivatives of the equation (\ref{KKT1}) with regard to $Q_p$ and $Q_r$ and then equating with zero, thus the solutions form as 
\begin{equation}\label{partialLagQp}
Q_p^{**}=\sqrt{\frac{2A_pD_p}{h_1+2\lambda_1 D_p C_3p_1+h_2pr+2\lambda_2 D_p C_3 p p_2 r}},
\end{equation}
and
\begin{equation}\label{partiaLaglQr}
Q_r^{**}=\sqrt{\frac{2\lambda C_2A_rD_r}{C_1C_2D_r(h_1+h_2)+2D_rh_2pr+\lambda C_1C_4}},
\end{equation}
where $\ds C_4=C_1C_2h_1+4h_2pr+2\lambda_2D_pC_3pp_2r+\frac{C_1C_2D_ph_2pr}{D_r}.$ \\ \\
Now we can look on the following four cases\par 
\noindent Case I: Putting $\lambda_1=\lambda_2=0$ in (\ref{partialLagQp}) and (\ref{partiaLaglQr}), we have (\ref{partialQp}) and (\ref{partialQr}). We can claim that (\ref{partialQp}) and (\ref{partialQr}) are  optimal solutions of (\ref{mathemodel}) if the solutions satisfy the constraints of the problem.\par 
\noindent Case II: If $\lambda_1\neq0$ and $\lambda_2=0$, then we obtain
\begin{equation*}
\lambda_1=\frac{2A_pD_p-Q_p^2(h_1+h_2pr)}{2Q_p^2D_pC_3p_1},
\end{equation*}
\begin{equation*}
Q_p^{**}=\frac{k_1}{p_1},
\end{equation*}
and $Q_r^{**}$ is as same as (\ref{partialQr}). If $\lambda_1 > 0$, then we can claim that these are global optimal solutions of (\ref{mathemodel}) if the solutions satisfy constraints of (\ref{mathemodel}).\par 
\noindent Case III: If $\lambda_1=0$ and $\lambda_2 \neq0$, then we make
\begin{equation}\label{case3a}
Q_p^{**}=\sqrt{\frac{2A_pD_p}{h_1+h_2pr+2\lambda_2 D_p C_3 p p_2 r}},
\end{equation}
and
\begin{equation}\label{case3b}
Q_r^{**}=\frac{D_r(k_2-rpQ_p^{**}p_2)}{rpD_pC_1p_2}.
\end{equation}
Plugging the solutions (\ref{case3a}) and (\ref{case3b}) into (\ref{KKT1}) evaluates $\lambda_2$. If $\lambda_2 > 0$ and the solution is feasible of the Problem~(\ref{mathemodel}), then global optimal solution is obtained.\par 
\noindent Case IV: If $\lambda_1\neq0$ and $\lambda_2 \neq0$, thereafter by complementary slackness conditions (\ref{KKT2}) and (\ref{KKT3}), we reach \begin{equation} \label{case4a}
Q_p^{**}=\frac{k_1}{p_1},
\end{equation}
and
\begin{equation}\label{case4b}
Q_r^{**}=\frac{D_r(k_2-rpQ_p^{**}p_2)}{rpD_pC_1p_2}.
\end{equation}
Combining the solutions (\ref{case4a}) and (\ref{case4b}) into the system (\ref{partialLagQp}) and (\ref{partiaLaglQr}) to assess $\lambda_1$ and $\lambda_2$. If $\lambda_1, \lambda_2 > 0$ and the solutions of the Problem~(\ref{mathemodel}) are feasible, then global optimal solution is obtained.

\section{Multiobjective Model and Solution Approach}\label{sec3}

Many environmental factors can arise in the reverse logistics inventory system. For example, greenhouse gas (GHG) emissions, substantial waste disposal pollution, and energy consumption can occur during the production of the products. It is required to control these environmental effects when one needs to optimize the average holding cost of the reverse logistics system. In our analysis, we initiate two more objectives such as minimization of greenhouse gas emission and energy used during the production and remanufacturing processes along with the objective of the total average holding cost described in Section~\ref{sec2}. \par 
\noindent Two variables are introduced to construct the mathematical model of the multiobjective reverse logistics problem, and these are
\begin{itemize}
	\item[]$Q_p$:  Procurement batch size\,,
	\item[]$Q_r$: Repair batch size\,.
\end{itemize}
Three objective functions for the above model are considered, one is introduced in~(\ref{costfunction}) as: 
\begin{itemize}
	\item [] $f_1(Q_p,Q_r)$= $\ds \frac{1}{C_3Q_p}(A_p+nA_r+h_1A_1+h_2A_2)$,
\end{itemize}
and other two functions listed in \cite{Bazan16} that illustrated as follows: \\
The greenhouse gas (GHG) emissions (ton per unit) are produced in the production process at a rate $\ds P=\frac{D_p}{M}$, where $\ds M=1-\frac{2A_pD_p}{h_1Q_p^2}>0$ introduced in \cite{Jaber2013a, Bogaschewsky1995} and therefore the second objective function is as follows:
\begin{itemize}
	\item[]$f_2(Q_p)$= $\ds \left(a_p\left(\frac{D_p}{M}\right)^2-b_p\left(\frac{D_p}{M}\right)+c_p\right)$, 
\end{itemize}
where $a_p$ an emissions function parameter (ton year$^2$/unit$^3$), $b_p$ an emissions function parameter (ton year/unit$^2$) and $c_p$ an emissions function parameter (ton year/unit) for manufacturing/production. It is here mentioned that, since repair/remanufacturing rate  $\lambda$ is constant so the second objective function $f_2$ does not have any effect by $Q_r$, therefore, it is only a function of $Q_p$.\par 	
\noindent And, the third objective function for energy consumed (KWh/year) in the production and remanufacturing process studied in \cite{Bazan2015} and \cite{Bogaschewsky1995}--\cite{Nolde2010} is defined as follows:
	\begin{itemize}
	\item[]$f_3(Q_p,Q_r)$= $\ds \frac{1}{T}\left[\left(\frac{MW_p}{D_p}+K_p \right)Q_p+\left(\frac{W_r}{\lambda}+K_r \right)nQ_r\right]$,
\end{itemize}
where $W_p$ and $K_p$ are respectively, idle power ($KW$/year) and  energy used (KWh/unit) for the manufacturing/production and also $W_r$ and $K_r$ are respectively, idle power ($KW$/year) and  energy used (KWh/unit) for the remanufacturing process. \par 
Thus, we construct the reverse logistics problem as a nonlinear multiobjective optimization problem as follows:
\begin{equation}\label{multimodel}
\begin{array}{rl} \min & \ [f_1, f_2, f_3]\\
\mbox{subject to the constraints} \\
&\ds  p_1Q_p \leq k_1, \\
& \ds p_2 \left(\frac{C_1Q_r}{D_r}+\frac{Q_p}{D_p}\right)rpD_p\leq k_2,\\
& \ds 1-\frac{2A_pD_p}{h_1Q_p^2} > 0.
\end{array}
\end{equation}
The solutions of (\ref{multimodel}) are called {\em efficient points} \cite{Yu1985},
or {\em Pareto points} \cite{Miettinen1999}. A less restrictive
concept of solution of (\ref{multimodel}) is the one of a {\em weak Pareto}
point. We use the following definitions of Pareto point and weak Pareto point.  
[see details in \cite{Chankong1983} and \cite{Miettinen1999}]. \par 
We intorduce the following standard notation to derive the definitions of efficient point and weak efficient point \cite{BurKayRiz2014}.  Let $u,v\in \mathbb{R}^{\ell}$ we
write that
\begin{equation*}\label{pareto-comp}
\begin{array}{c}
u\leqq v \text{ if and only if } u_i\leq v_i
\,\forall\,i=1,\ldots,\ell\,\text{ and } \\
\,\exists\,j\text{ such that } u_j< v_j\,;\\
u<v \text{ if and only if } u_i< v_i
\,\forall\,i=1,\ldots,\ell\,.
\end{array}
\end{equation*}
 Let $X$ be a feasible set of Problem~(\ref{multimodel}). A point $\bar{x} \in X$ is said to be
		\textit{Pareto} for Problem (\ref{multimodel}) iff there is no $x \in X$, $x \neq \bar{x}$ such that $\ f(x) \leqq f(\bar{x})$.
		 Moreover, a point $\bar{x} \in X$ is said to be		\textit{weak Pareto} for Problem (\ref{multimodel}) iff there is no $x \in X$
		such that $\ f(x) < f(\bar{x})$. 
\subsection{Solution Approach of the Proposed Multiobjective Problem}\label{sec3a}
In this section, we recall the scalarization approach introduced in \cite{BurKayRiz2017} named the Objective-Constraint Approach to solve the proposed Problem (\ref{multimodel}). The method seems to perform efficiently to construct boundary and interior of the Pareto front when the front is disconnected.  For more details on these techniques, see \cite{BurKayRiz2017} and references therein.
\par
\textbf{The Objective-Constraint Approach}: For $\hat{x} \in X$, define the weight
\begin{equation}\label{weiassum}
\hspace*{-5mm} w_i:=\frac{1/f_i(\hat{x})}{\ds \sum_{j=1}^{\ell}1/f_j(\hat{x})}.
\end{equation}
Then $w \in W$ and $w_kf_k(\hat{x})=w_if_i(\hat{x}), \; \forall i=1,2,..., \ell, \;\; i \neq k.$
The associated scalar problem is defined as
\begin{equation*} \label{obj_con}
\hspace*{-5mm} \mbox{($P_{\hat{x}}^k$)} \hspace{20mm}\;\; \ \left\{\begin{array}{rl} \ds\min_{x\in  X} & \
\ w_{k}f_{k}(x),
\\[4mm]
\mbox{subject to} & \ \ w_if_i(x) \leq w_kf_k(\hat{x}), \;\;\mbox{i=1,2,...,$\ell$}, \;\; i \neq k.
\end{array}
\right.
\end{equation*}
Note that if $\bar{x}$ is an efficient solution of Problem (\ref{multimodel}) and $w$ is as in (\ref{weiassum}), then $\bar{x}$ is the solution of ($P_{\bar{x}}^k$) for all $k$. However, the converse is not true. For weak efficient point the above property holds for both necessary and sufficient conditions.
\section{ Numerical Experiments, Results and Analysis}\label{sec4}
\begin{figure}[hbt!]
	\hspace{-10mm}
	\begin{minipage}{90mm}
		\begin{center}
			\includegraphics[width=85mm]{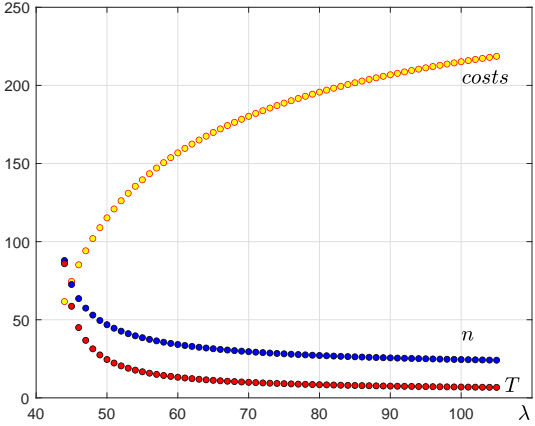} \\
			{\scriptsize (a) Holding costs, Repair cycle $n$, and \\ \hspace*{-8mm} Time cycle $T$  at $D_r=43$.}
		\end{center}
	\end{minipage}
	\hspace{0mm}
	\begin{minipage}{90mm}
		\begin{center}
			\hspace{-20mm}
			\includegraphics[width=80mm]{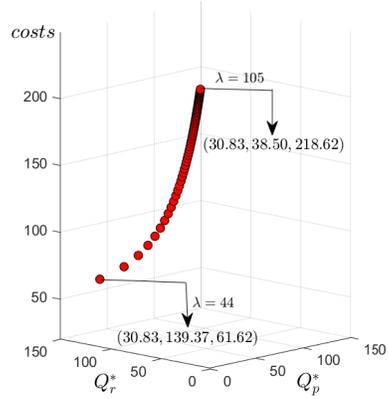} \\
			{\scriptsize (b) Optimal solution $(Q_p^*,Q_r^*, costs)$\\ \hspace*{5mm} at $D_r=43$ and for $44<\lambda<105$.}
		\end{center}
	\end{minipage}
	\\[2mm]
	\hspace*{-1cm}
	\begin{minipage}{90mm}
		\begin{center}
			\hspace*{0cm}
			\includegraphics[width=80mm]{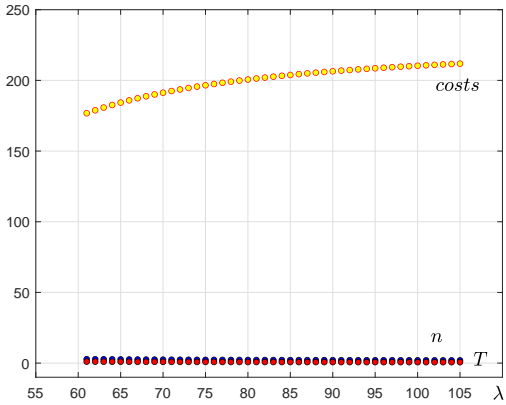} \\
			{\scriptsize (c) Holding costs, Repair cycle $n$, and \\ \hspace*{-8mm}Time cycle $T$ at $D_r=60$.}
		\end{center}
	\end{minipage}
	\hspace{-0cm}
	\begin{minipage}{90mm}
		\begin{center}
			\hspace*{-2cm}
			\includegraphics[width=85mm]{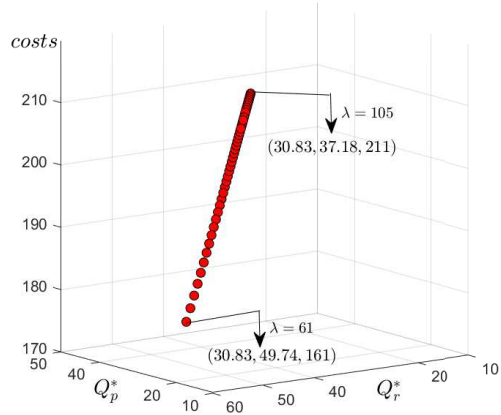} \\
			{\scriptsize (d) Optimal solution $(Q_p^*,Q_r^*, costs)$\\ \hspace*{5mm} at $D_r=60$ and for $61<\lambda<105$.}
		\end{center}
	\end{minipage}
	\caption{Optimal solutions with average holding costs, repair and time cycles approximations over the selection of repair rates for Example--\ref{expuncon}. Figures (a)(c) are presented relations among costs, repair cycle $n$ and time cycle $T$ when $D_r=43$ and $D_r=60$, respectively, whereas optimum solutions and associated average holding costs $(Q_p^*,Q_r^*, costs)$ are depicted in Figures (b) (d) that are approximated  for  $D_r=43$ and $D_r=60$, over the arbitrary interval of $\lambda$.}
	\label{cost_n_T_uncon}
\end{figure}

\begin{figure}[hbt!]
		\begin{center}
			\hspace*{-15mm}
			\psfrag{f1}{$f_1$}
			\psfrag{f2}{$f_2$}
			\includegraphics[width=120mm]{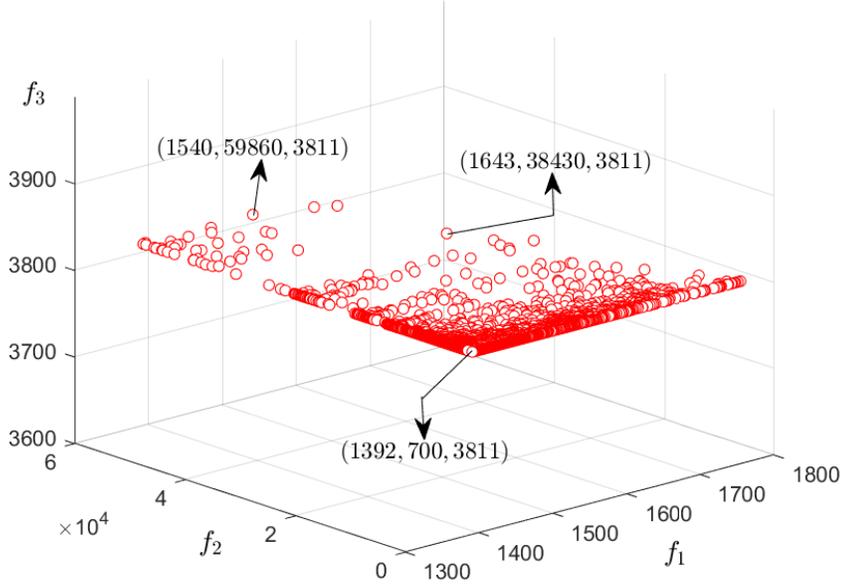} \\[3mm]
		\end{center}
	\caption{ \small{ The Pareto points are obtained for Example~\ref{multexam}}.}
	\label{figmulti}
\end{figure}

In this section, we test our proposed models that are stated in (\ref{costfunction}), (\ref{mathemodel}) and (\ref{multimodel}), and the solution formulas introduced in (\ref{partialQp})--(\ref{partialQr}) and (\ref{partialLagQp})--(\ref{partiaLaglQr})  for the range of input parameter settings in the tyre industry. We first introduce
 Example-\ref{expuncon} to test the  average holding cost function (\ref{costfunction}) which depends on a set of parameters. In this problem, we consider the unlimited floor spaces in supply and repair depots, and calculate average holding costs of procurement and repair items in per unit time of cycle $T$. Example-\ref{exm_con} is  commenced for certain situations where the restriction of floor spaces in supply and repair depots need to be considered. Therefore, we test constraint-model (\ref{mathemodel}) for the same input parameters that are used in Example-\ref{expuncon}. We also introduce a more challenging problem for the model (\ref{multimodel}) where multiple objective functions are considered. As far we know, not enough models and algorithms are presented in the literature of the reverse logistics systems where more objective functions are taken, and solutions are reported. By considering this fact, we presented here Problem-\ref{multexam} to analyze the proposed multiobjective model (\ref{multimodel}), and approximate the solution set of the Pareto front in Figure~\ref{figmulti}. 
\begin{example}\label{expuncon}
	 Suppose that in the tyre industry, tyre procurement and repair setup costs are $A_p=\$ \;10$ and $A_r=\$ \;30$ for batch sizes $Q_p$ and $Q_r$, respectively. The demand of the new and repaired tyres are $D_p=100$ and $D_r=43$, respectively, over cycle $T$. We also assume that the collection percentage of available returns of used items is $p=0.6$, recovery rate is $r=0.7$ and  repair rate is $\lambda$.  It is here mention that $D_r$ is smaller than $\lambda$, and they both are always larger than $rpD_p$ according to our proposed model and as above $rpD_p=42$. Let us also assume that holding costs  per unit per unit of time for supply and repair depots are $h_1=\$\;1.6$ and $h_2=\$\;1.2$, respectively. We intend to find $Q_p$ and $Q_r$ so that the average holding cost would be the minimum over per unit time of cycle $T$. Therefore, we aim to minimize the model (\ref{costfunction}) under the above parameter settings.
\end{example}
	Now we employ the obtained solution formulas (\ref{partialQp}) and (\ref{partialQr}) for the Example-\ref{expuncon}. According to the assumptions of the model (\ref{costfunction}), in our experiment, we first set $\lambda=45$ which is always larger than $rpD_p$ and $D_r$.   As a result, when $D_r=43$ and $\lambda=45$ we attain the  optimal procurement and repaired batches are $Q_p^*=30.83$ and $Q_r^*=115.10$ with the optimal  average holding cost $74.61$, which are dipicted in Table~\ref{table:exp-2a} . In our analysis, we observe that, the optimal cost increases when $\lambda$ increases ($\lambda >44$),  which is shown in Figure \ref{cost_n_T_uncon}(a). With the same $D_r=43$, when $\lambda$ increases gradually to $105$, the optimal solution and 'average holding çosts' are changes to $(Q_p^*,Q_r^*, holding    
	costs)=(30.83, 40.83, 218.63)$, these changes are presented in Figure \ref{cost_n_T_uncon}(b). Besides when $\lambda=45$ and $D_r=43$, the repair cycles is $72$ (consider the integer number) and the time cycle is $58.62$ (see, in Table~\ref{table:exp-2a}), however, when $\lambda$ increases these two cycles are decreases (see, Figure \ref{cost_n_T_uncon}(a)). The above analysis indicates that, when repair rate increases for fixed $D_r$ then the inventory of repaired items in supply depot increases, and therefore the average holding costs of the items are also increases. Here, decision makers have the option of choosing an appropriate $\lambda$ to obtain average holding costs and cycle loops, which may require a trade-off with the optimum average holding costs. Similar results obtain when $D_r=60$, this can be seen in Figure \ref{cost_n_T_uncon}(c)(d).
  \par
	\begin{table}
		\caption{\small{\textit{Example-\ref{expuncon} -- Numerical performance of model \eqref{costfunction}, and formulas (\ref{partialQp}) and (\ref{partialQr}). }}}
		\footnotesize
		\vskip 1.5em
		\centering
		\begin{tabular}{|c| c| c| c| c| c| c| c| c|}
			\hline
			Repair &     Repair 		&  Demand  for  	& Procurement & Repaired& Holding& Repair &Time& CPU   \\
			portions &     rates 		&   repaired items  	& batch sizes & batch sizes& costs& cycles &cycles & time  \\
			$rpD_p$ &     $\lambda$ 		&  $D_r$   	& $Q_p^*$ & $Q_r^*$& ${\tiny f_1(Q_p^*,Q_r^*)}$ & $n$ &$T$&  [sec] \\ [0.5ex]
			\hline
			42&45 &43 & 30.83  &	115.10	&74.61 		& 72.56 & 58.62 & 1.11  \\ [.5ex]\hline
			42	&60 &43 & 30.83 &54.53 		&156.81 &	34.15	& 13.20 & 0.73  \\ [.5ex]\hline
			42	& 75& 43 &30.83 & 44.92 		&188.68 &	28.17	& 9.07 		& 0.88\    \\ [.5ex]
			\hline
			42	& 90& 43 &30.83 & 40.83 		&206.80 &	25.75	& 7.52 		& 0.63\    \\ [.5ex]
			\hline
			42	&105 &43 &30.83 &38.51  		&218.63 	& 24.10	& 6.70 & 0.61  \\ [.5ex]\hline
		\end{tabular}
		\label{table:exp-2a}
	\end{table}

\begin{figure}[hbt!]
	\begin{minipage}{85mm}
		\begin{center}
			\includegraphics[width=85mm]{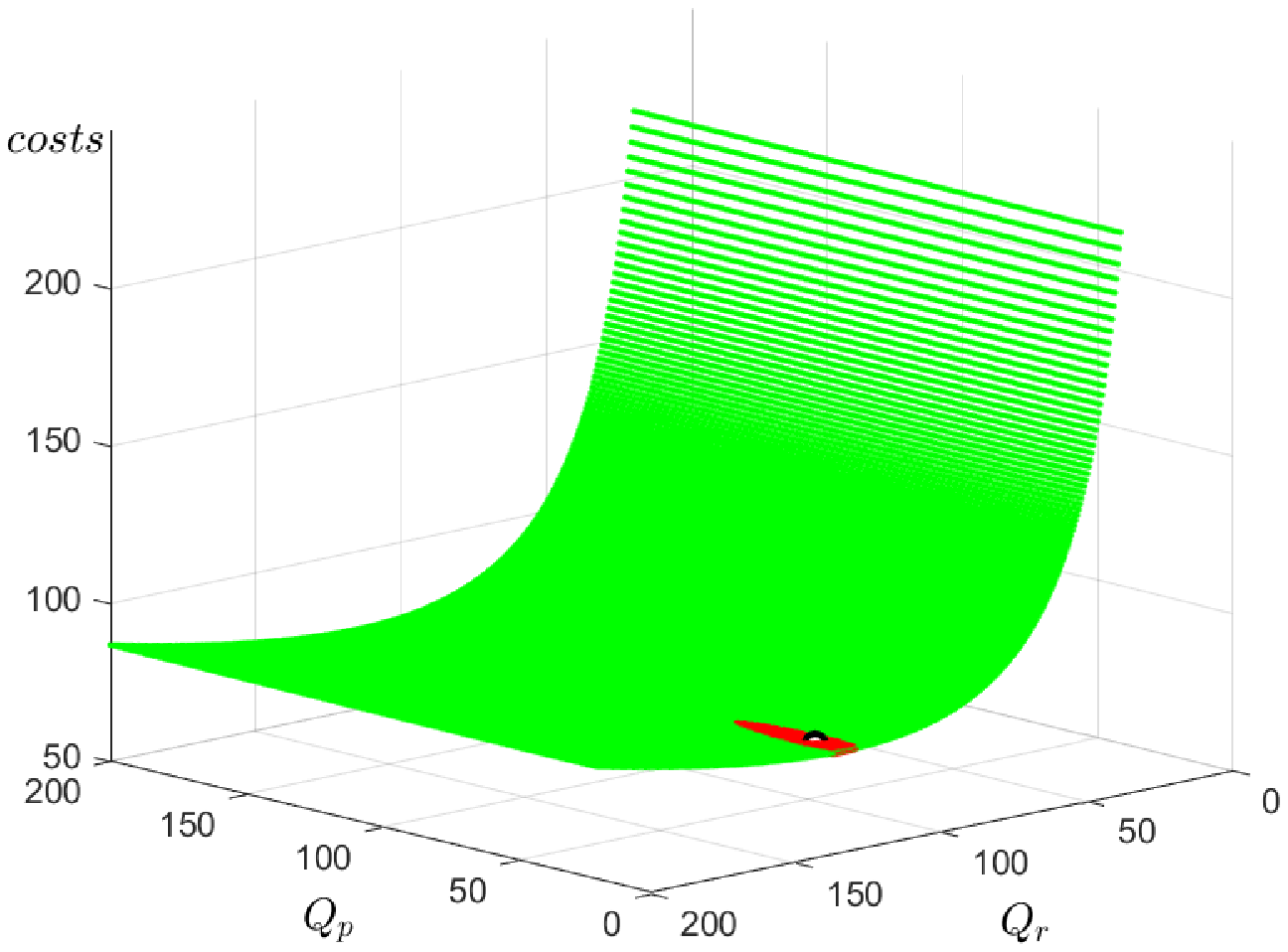} \\
			{\small (a) Solution position on the surface \\ when $\lambda = 45$ and $D_r=43$}.\\
		\end{center}
	\end{minipage}
	\begin{minipage}{85mm}
		\begin{center}
			\hspace*{-15mm}
			\psfrag{f1}{$f_1$}
			\psfrag{f2}{$f_2$}
			\includegraphics[width=85mm]{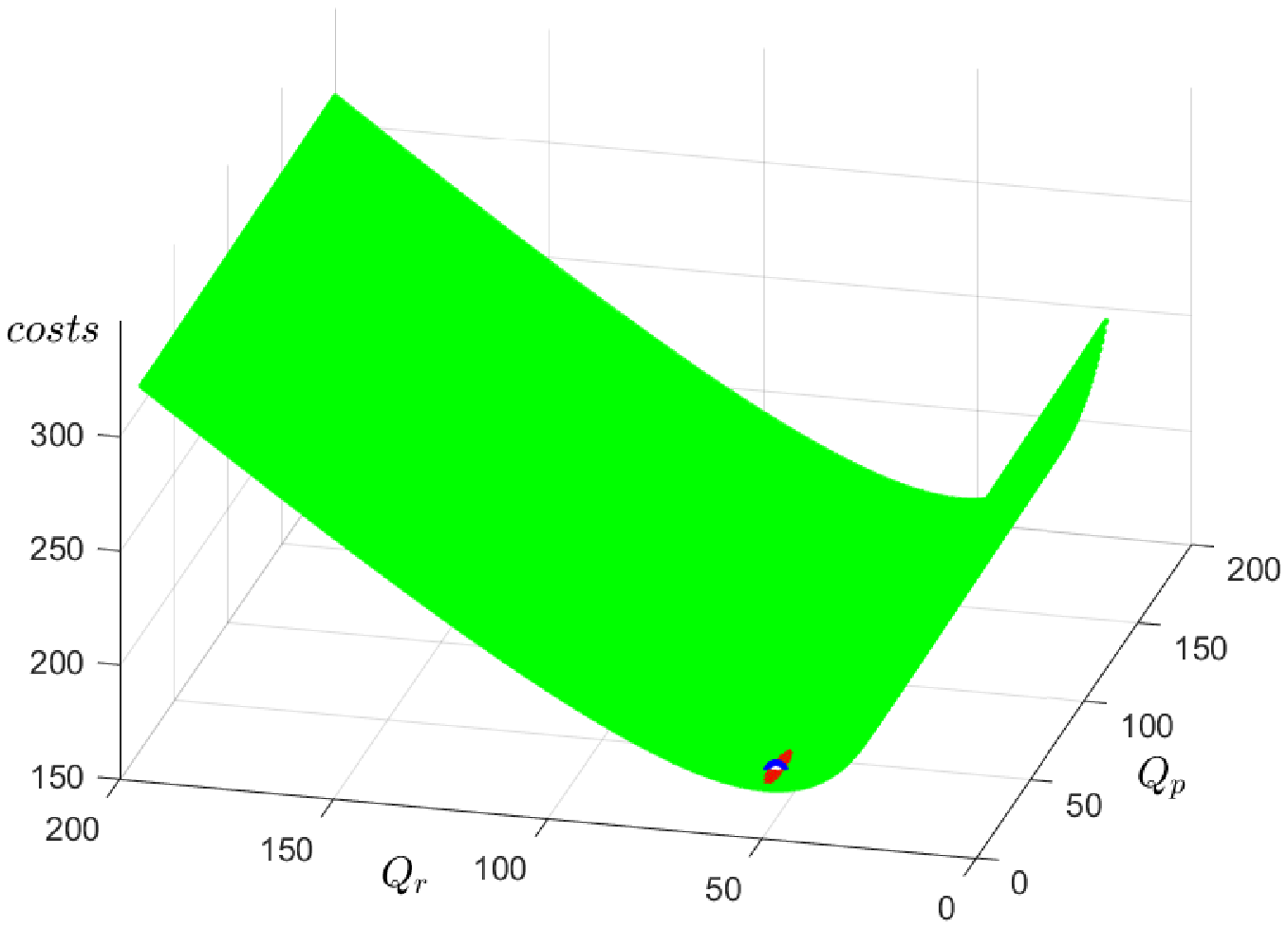} \\
			{\small(b) Solution position on the surface \\ when $\lambda = 60$ and $D_r=43$}.
		\end{center}
	\end{minipage}
\\[2mm]
\hspace*{-1cm}
\begin{minipage}{100mm}
	\begin{center}
		\hspace*{0cm}
		\includegraphics[width=90mm]{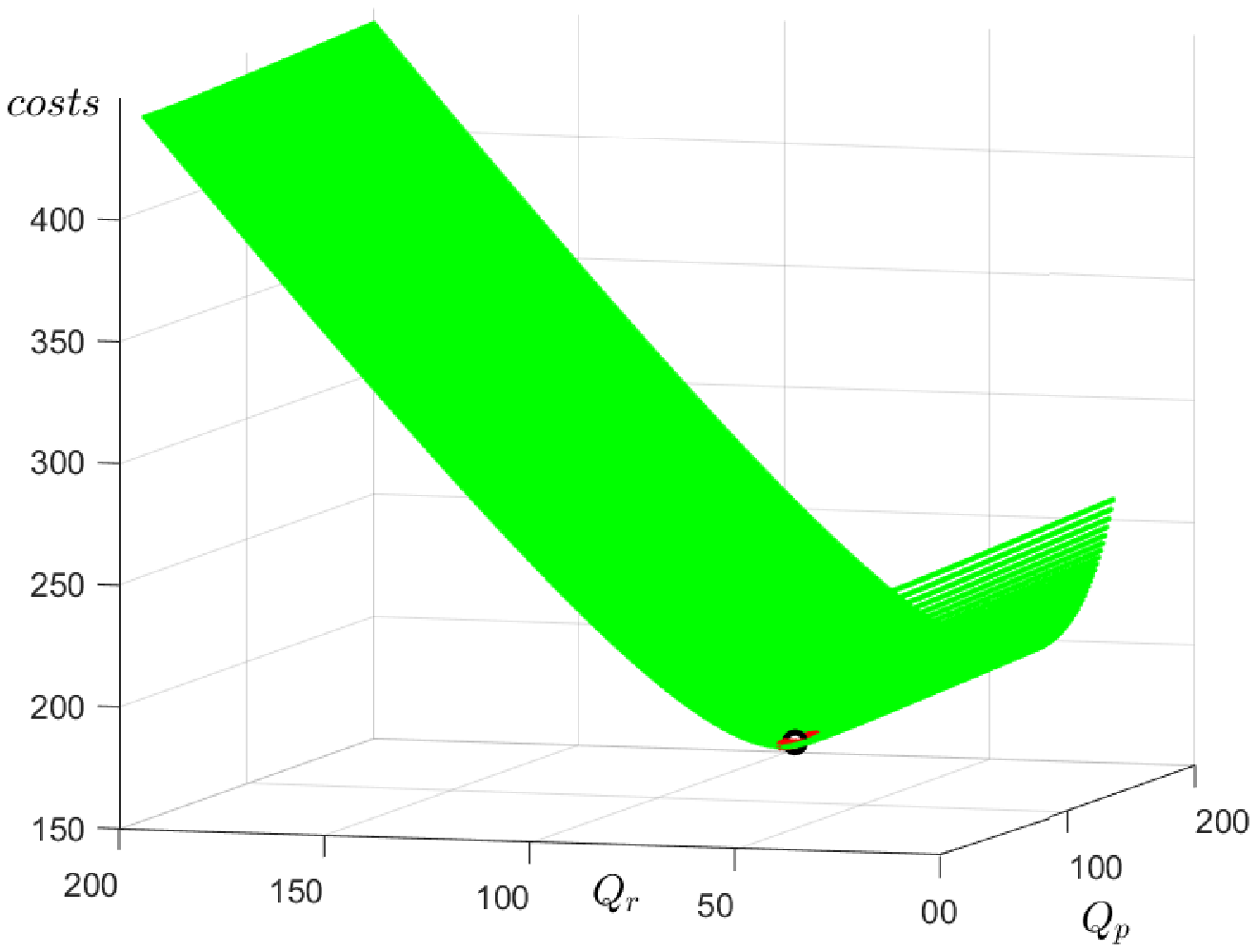} \\
		{\small(c) Solution position on the surface \\ when $\lambda = 75$ and $D_r=43$}.
	\end{center}
\end{minipage}
\hspace{-0cm}
	\caption{ \small{(Green) surface produce by (\ref{costfunction}). (Red) circles represent the solution points of Example-\ref{expuncon} were obtained by solvers and (Blue) big circle indicates a solution obtained by  (\ref{partialQp}) and (\ref{partialQr}).}}
	\label{figunconR}
\end{figure}

In our experiments, we employ other kinds of techniques to approximate the solution of the Example-\ref{expuncon}, for instance, Example-\ref{expuncon} is solved using a range of nonlinear solvers such as fmincon with sequential quadratic programming, Ipopt \cite{Wachter2006}, and SCIP \cite{Achterberg2009}. We observe that the approximated results vary for solvers and initial guesses. One can be interested to see the whole set of solutions of Example-\ref{expuncon}. Therefore, we utilize a search base algorithm, like Brute Force algorithm, to approximate the solution set. In Figure~\ref{figunconR}, the straight line formed by the small circles (red) is the approximation of the solutions of Example-\ref{expuncon} obtained by search based algorithm. \par 
\begin{figure}[hbt!]
	\begin{minipage}{85mm}
		\begin{center}
			\hspace*{-15mm}
			\psfrag{f1}{$f_1$}
			\psfrag{f2}{$f_2$}
			\includegraphics[width=80mm]{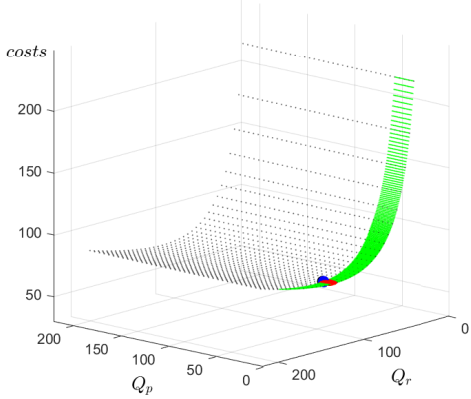} \\
			\small {(a) Solution position on the surface \\ when $\lambda = 45$ and $D_r=43$}.
		\end{center}
	\end{minipage}
	\begin{minipage}{85mm}
		\begin{center}
			\hspace*{-15mm}
			\psfrag{f1}{$f_1$}
			\psfrag{f2}{$f_2$}
			\includegraphics[width=80mm]{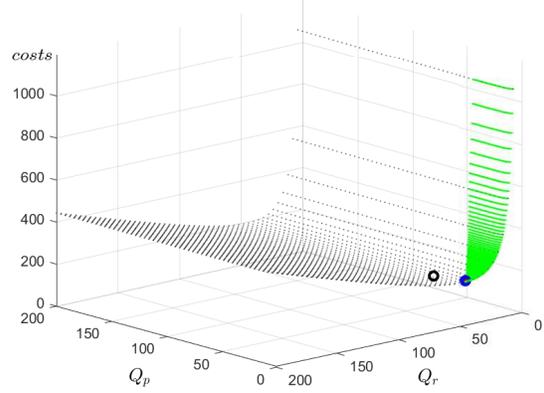} \\
			\small{(b) Solution position on the surface \\ when $\lambda = 75$ and $D_r=43$}.
		\end{center}
	\end{minipage}
	\caption{ \small{(Green) surface is the feasible space of Example-\ref{exm_con} whereas black surface obtained from Example-\ref{expuncon} is discarded here because of the floor restrictions. (Red) circles represent the solution points were obtained by solvers and (Blue) big circle indicates a solution obtained by  (\ref{partialLagQp}) and (\ref{partiaLaglQr}).}}
	\label{figcon}
\end{figure}
\begin{figure}
	\hspace{-10mm}
	\begin{minipage}{90mm}
		\begin{center}
			\includegraphics[width=85mm]{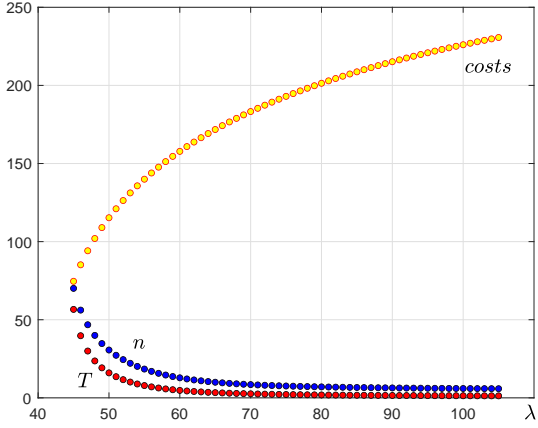} \\
			{\small (a) Holding costs, Repair cycly $n$, and\\ Time cycle $T$.}
		\end{center}
	\end{minipage}
	\hspace{0mm}
	\begin{minipage}{90mm}
		\begin{center}
			\hspace{-20mm}
			\includegraphics[width=85mm]{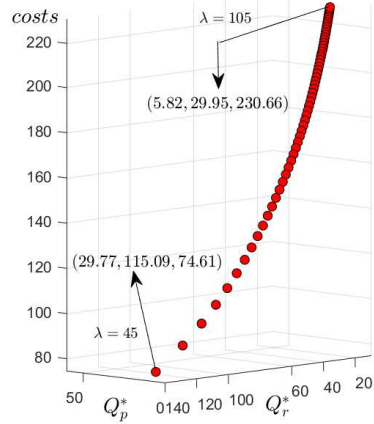} \\
			{\small (b) Optimal solution $(Q_p^*,Q_r^*, costs)$.}
		\end{center}
	\end{minipage}
	\caption{Figure (a) is presented optimal average holding costs, repair cycle $n$ and time cycle $T$ over the interval for repair rate $44<\lambda<105$  when $D_r=43$ for Example--\ref{exm_con}, whereas optimum solutions $(Q_p^*,Q_r^*, costs)$ are depicted in Figure (b) for the same $\lambda$'s and $D_r$.}
	\label{conDr}
\end{figure}
Now we introduce constrained Example-\ref{exm_con} and test the model (\ref{mathemodel}) and its proposed solution (\ref{partialLagQp}) and (\ref{partiaLaglQr}). 
\begin{example} \label{exm_con}	Consider the input parameters that are used in Example-\ref{expuncon}.
Moreover, we assume that the amount of square feets required to each of the item in supply and repair depots are $p_1=0.5$ and $p_2=0.5$, respectively. We also assume that the availability of maximum number of floor spaces for supply and repair depots are, $k_1=20$, $k_2=10$, in feets, respectively. We intend to find $Q_p$ and $Q_r$ so that the total  average holding cost would be the minimum per unit time of cycle $T$. Therefore, we are interested to minimize the model (\ref{mathemodel}) under the above parameter settings.
\end{example}
In the first instance, if we consider $D_r=43$, $\lambda=45$ and $rpD_p=42$, then formulas (\ref{partialLagQp})--(\ref{partiaLaglQr}) give the solution of Example-4.2, and thus Case I of the model (\ref{mathemodel})  provides the optimal procurement and repaired batches which are $Q_p^*=29.77$ and $Q_r^*=115.09$, and the optimal  average holding  cost is $74.61$ (see Table~\ref{table:exp-2b}). Whereas we obtain an alternative solution $(30.5,107.5)$ with the approximate optimal  average holding cost  $74.78$  from Case II. Lagrange multiplier conditions are not satisfied for  Cases III and IV therefore, the solutions are discarded. Furthermore, for fixed $D_r=43$ and the optimal batch sizes $Q_p^*$ and $Q_r^*$, we also calculate repair and time cycles which are  $n=70$ (taken only integer part) and $T=56.61$. When $\lambda$  increases within an arbitrary interval $44\leq \lambda \leq 105$ then optimal average holding costs increases, however two cycles $n$ and $T$ are decreases which are depicted in Figures~\ref{conDr}(a)(b). Therefore, one can cautiously choose the repair rate compare to the demand rate $D_r$ of repaired items, this might be required to trade-off among  average holding costs, repair and time cycles. On the other hand, for  $D_r=43$ and  $\lambda$ that are chosen randomly from $\ds \left\{60, 75, 90, 105\right\}$, the optimal solutions and associated average holding costs are changed that are demonstrated in Table~\ref{table:exp-2b}.  The above analysis indicates that when the repair rate is  gradually increased that are far larger than $D_r$, then the optimal average holding costs would be increased, whereas both repair and time cycles significantly are decreased (see Figure \ref{conDr}(a)).  These experiments have also been conducted by using solvers such as fmincon with sequential quadratic programming and Ipopt \cite{Wachter2006}, and we obtain the same set of approximations of the solution. Brute Force algorithm is also used to obtain the whole set of solutions. The solution set obtained by the experiments  with Table~\ref{table:exp-2b} is depicted in Figure~\ref{figcon} as (red) small circles. Note that, in these input parameter settings, the obtained feasible space and optimal solutions of constrained model (\ref{mathemodel}) are not the same as Example-\ref{expuncon}, because floor restrictions are applied in model (\ref{mathemodel}).
\begin{table}
	\caption{\small{\textit{Example-\ref{exm_con} -- Numerical performance of model \eqref{mathemodel}, and formulas (\ref{partialLagQp}) and (\ref{partiaLaglQr}).}}}
	\footnotesize
	\vskip 1.5em
	\centering
	\begin{tabular}{|c| c| c| c| c| c| c| c| c|}
		\hline
		Repair &     Repair 		&  Demand  for  	& Procurement & Repaired& Holding& Repair &Time& CPU   \\
		portions &     rates 		&   repair items  	& batch sizes & batch sizes& costs& cycles &cycles & time  \\
		$rpD_p$ &     $\lambda$ 		&  $D_r$   	& $Q_p^*$ & $Q_r^*$& $f_1(Q_p^*,Q_r^*)$ & $n$ &$T$&  [sec] \\ [0.5ex]
		\hline
	42&45 &43 & 29.77  &	115.09	&74.61 		& 70.08 & 56.61 & 1.59  \\ [.5ex]\hline
		42	&60 &43 & 11.13 &52.29 		&157.78 &	12.82	& 4.77 & 1.01  \\ [.5ex]\hline
		42	& 75& 43 &7.28 & 39.42 		&193 &	7.58	& 2.14 		& 0.65\    \\ [.5ex]
		\hline
		42	& 90& 43 &6.26 & 33.35 		&215.15 &	6.35	& 1.53 		& 1.1\    \\ [.5ex]
		\hline
		42	&105 &43 &5.82 &30  		&230.7 	& 5.85	& 1.27 & 0.914  \\ [.5ex]\hline
	\end{tabular}
	\label{table:exp-2b}
\end{table}


 In addition, there are requirements to minimize the environmental pollution that arises during the production and remanufacturing process in tyre industry. The environmental impacts from tyre production include greenhouse gas emissions, energy consumption, dust emission and solvent emission etc.  We consider two additional objective  functions with Example-\ref{exm_con} to minimize greenhouse gas emission and energy consumption  during the production and remanufacturing process, which is presented here as Example-\ref{multexam}. \par

In Example-\ref{multexam}, three objective functions are considered which are  average holding cost function $f_1(Q_p;Q_r)$, greenhouse gas emissions function $f_2(Q_p)$ and energy consumption function $f_3(Q_p;Q_r)$. Three objective functions with the constraints set make the Example (\ref{multexam}) of model (\ref{multimodel}) hard to solve. Moreover, the problem has disconnected Pareto front, therefore, we need to choose a scalarization method carefully to solve the problem. We choose ($P_{\hat{x}}^k$) to solve the problem as the method is efficient to approximate the Pareto front when the Pareto front is disconnected. 

\begin{example}\label{multexam}
	We make changes in the parameter settings that are used in Examples \ref{expuncon} and \ref{exm_con}  as the available differtial solvers have difficulties to approximate solutions for certain settings. In this Example we assume that $r=0.7$, $p=0.6$, $D_p=1000$, $D_r=422$; $\lambda=450$, $A_p=50$, $A_r=100$, $h_1=20$, $h_2=10$, $p_1=1$, $p_2=1$, $k_1=2000$, $k_2=2000$. We set emission function parameters as $a_p=0.00000003$ (ton year$^2$/unit$^3$), $b_p=0.0014$ (ton year/unit$^2$) and $c_p=1.4$ (ton year/unit) for production process. Let, the idle powers $W_p=120$ ($KW/year$)  and $W_r=80$ ($KW/year$) for production and remanufacturing process, respectively. Moreover, we set energy $K_p=5.5$ $(KWh/unit)$, and $K_r=2.5$ $(KWh/unit)$  are used  in the production and remanufacturing process, respectively. 
 We aim to approximate the Pareto points of the multiobjective optimization model (\ref{multimodel}) under the above parameter settings.
\end{example}

We now introduce the following Algorithm to solve Example-\ref{multexam}.  

\subsection*{Algorithm}

We use the objective-constraint approach ($P_{\hat{x}}^k$) in the algorithm.  In Step 2 of Algorithm, each objective function is minimized subject to the original constraints of the Example-(\ref{multexam}). These individual optimum point are used in Step 3 to form weighted grids. Each grid point corresponds to a weight vector in $\mathbb{R}^3$. In Step 4, three sub-problems are solved at each grid point to generate Pareto points. \par
In Step~4(b), we calculate the efficient and weak efficient points using the fact that if $\bar{x}_1=\bar{x}_2=\bar{x}_3=:\bar{x}\; \mbox{(say)}$ holds, then the solution $\bar{x}$ is an efficient. Here,  $\bar{x}_k$ are the solutions of ($P_{\hat{x}}^k$) for $k=1,2,3$. On the other hand, if $\bar{x}_1=\bar{x}_2=\bar{x}_3$ does not hold, then any dominated point is removed from the set \{$\bar{x}_1$, $\bar{x}_2$ and $\bar{x}_3$\} (see Step 4(b)), and these are all weak efficient points. The latter case is typically encountered when the Pareto front and/or the domain is disconnected. Therefore, this algorithm is efficient in finding Pareto points even when the feasible set is discrete or disconnected.

\begin{description}
	\item[Step $\mathbf{1}$] { \textbf{(Input)}} \\ Set all parameters that stated in Example-\ref{multexam}.
	\item[Step  $\mathbf{2}$] {\textbf{(Determine the individual minima)}}\\
	Solve Problem $\ds \min f_i, \;\; i=1,2,3$, subject to the constraints of Example-\ref{multexam}  that give the solutions $\left(\bar{Q}_{p_i}, \bar{Q}_{r_i}\right)$, for $i=1,2,3$, respectively.	
	\item[Step $\mathbf{3}$] { \textbf{(Generate weighted parameters)}} \\
	Weighted parameters are generated in this step. We generate weighted grids as introduced in \cite[Step 3 of Algorithm 3]{BurKayRiz2019}
	
	\item[Step $\mathbf{4}$] Choose $w=(w_1,w_2,1-w_1-w_2)$, which generated from Step 3.
	\begin{description}
		\item[(a)] Find $\hat{x}_k=\left(Q_{p_k}, Q_{r_k}\right)$ that solves Problems ($P_{\hat{x}}^k$), $k=1,2,3.$
		\item[(b)] Determine weak efficient points : 
				\begin{description}
			\item[(i)] If $\bar{x}_1=\bar{x}_2=\bar{x}_3$, then set $\bar{x}=\bar{x}_1$ (an efficient point)\\
			and  Record the points.
			\item [(ii)] If $\bar{x}_1=\bar{x}_2=\bar{x}_3$  does not hold, then, any dominated point is discarded by comparing these $3$ solutions. \\
			Record non dominated points.	
		\end{description}
	\end{description}
	\item[Step $\mathbf{5}$] (Output)\\
	All recorded points are Pareto point of Problem-\ref{multexam}.\\
\end{description}
We write code in MATLAB. We test the range of solvers in Steps 3 and 4(a), these include differential and non-differential solvers such as fmincon with sequential quadratic programming algorithm, SCIP \cite{Achterberg2009} and SolvOpt. We   take the input parameter values that are listed in the problem statement and implement the proposed algorithm. The parameter setting plays an important role in Example- \ref{multexam}. We choose input parameters in our problem randomly, and for the given parameter setting, we provide $1325$ weight-grids ($w$) into the algorithm. As a result, algorithm approximates $2527$ Pareto points which are depicted in Figure~\ref{figmulti}.  The elapsed CPU time was about $5$ minutes. Note that the computations have been performed on a DELL Inspiron 15 7000 laptop with 8 GB RAM and core i7 at 4.6GHz.
\section{Conclusion}
We developed a new material flow model in the reverse logistics system to measure holding cost under assumptions that the demands of new and repaired items are different and deterministic. We established mathematical expressions to compute the holding cost and provided formulae to optimize the holding cost. We extended the proposed model to the case where constraints were imposed, and the solution approaches were established. Extensive
computational experiments were conducted to demonstrate the efficiency of the proposed
models. MATLAB was used to perform the tests, and a wide range of solvers have been utilized to solve the complex nature problem. Moreover, we extended the proposed single objective problem into a three-objective problem that simultaneously optimizes holding
cost, greenhouse gas emissions, and energy consumptions during the production and remanufacturing processes. Well-known scalarization method employed to solve this three-objective optimization problem. We successfully approximated
the Pareto front of this problem, and the computational time has been reported.
\bigskip

\textbf{\large {Acknowledgments:}} 
The first author is supported by NST (National Science and Technology) fellowship, reference no. 120005100-3821117 in the session: 2019-2020, through the Ministry of National Science and Technology, Bangladesh. This financial support of NST fellowship is gratefully acknowledged.

\end{document}